\documentclass[11pt]{article}
\usepackage{latexsym}

\title{Riesz endomorphisms of Banach algebras}

\author{Joel F. Feinstein and Herbert Kamowitz}

\begin{document}
\newcommand{\dis}{\displaystyle}
\newcommand{\pis}{\mbox{$\overline{\phi(b)}$}}
\newcommand{\QED}{\hfill$\Box$\vskip.3cm}
\def\unorm#1#2{{\|#1\|_{{#2}}}}
\def\norm#1{{\|#1\|}}
\def\rd{{\rm d}}

\maketitle

\begin{abstract}
   Let $B$ be a unital commutative semi-simple Banach algebra.
 We study endomorphisms of $B$
which are simultaneously Riesz operators. 
Clearly compact and  power compact endomorphisms are
Riesz. Several general theorems about Riesz endomorphisms are proved,
and these results are then applied to the question of when 
Riesz endomorphisms of certain algebras are necessarily power compact.
\end{abstract}

\begin{center}
Introduction
\end{center}

     Let $B$ be a unital commutative semi-simple Banach algebra. An
endomorphism $T:B \rightarrow B$ is a linear operator which also preserves
multiplication. An endomorphism $T$ is called unital if $T1=1$.
Since the Banach algebra is assumed to be semi-simple, it follows from
very early theory that $T$ is necessarily bounded. It is interesting to
see what can be deduced if various operator theoretic properties are imposed
on the endomorphism. In this note we consider endomorphisms which are
also Riesz operators. We will start out by proving two theorems which are
valid for all unital commutative semi-simple Banach algebras, and then
apply them to specific examples.

  If $B$ is a Banach space, a Riesz operator $T:B \rightarrow B$ is a
bounded linear map satisfying properties (a) - (e) \cite{heus}.

(a) For each non-zero $\lambda$, $\lambda - T$ is open.

(b) For each non-zero $\lambda$, $(\lambda - T)(B)$ is closed.

    If $S$ is a linear operator $B \rightarrow B$, let $N(S)=\{f:Sf=0\}$,
and $R(S)=S(B)$ the range of $S$.

(c) For each non-zero $\lambda$, dim$N(\lambda - T)$ and codim$R(\lambda - T)$
are finite and equal.

(d) For each non-zero $\lambda$, the lengths of the null chain
$N(\lambda - T) \subseteq N((\lambda - T)^2) \subseteq \cdots,$
and the image chain $R(\lambda - T) \supseteq R((\lambda - T)^2) \cdots$
are finite and equal.

(e) The non-zero spectrum of $T$ consists of eigenvalues, and if there
are infinitely many of them, they form a sequence approaching $0$.

   Proposition 52.2 in \cite{heus} contains the following useful
characterization of Riesz operators.

{\bf Theorem:} Let $B$ be a Banach space and $T$ a bounded linear
operator from $B$ to $B$. Then $T$ is a Riesz operator if, and only
if,

$\dis \lim[\inf\{\|T^n - K\|: K$ is a compact operator $B \rightarrow B\}]^{1/n}=0.$

   The limit on the left hand side of the last equation is called the
essential spectral radius of $T$. Clearly compact operators and quasinilpotent
operators are Riesz operators. The sum of a compact operator and a
quasinilpotent operator is Riesz. A bounded linear operator $T$ is power
compact if for some positive integer $N$, $T^N$ is compact. Power
compact operators are Riesz. One of the questions  considered is for which
algebras is every Riesz endomorphism necessarily power compact.

\begin{center}
Part 1
\end{center}

    If $B$ is a unital commutative simi-simple Banach algebra with maximal
ideal space $X$, and $T$ is a unital endomorphism of $B$, then it is
well known that there exists a w*-continuous selfmap $\phi$ of $X$
such that for all $f \in B$ and $x \in X$,
$\widehat{Tf}(x) = \hat{f}(\phi(x))$.  In this case we  say that
$T$ is induced by $\phi$ or that $\phi$ induces $T$.

    It was shown in \cite{comp1} that if $\phi$ induces a compact
endomorphism of the unital commutative simi-simple Banach algebra $B$,
and the maximal ideal space $X$ of $B$ is connected, then
$\bigcap_{n=0}^\infty \phi_n(X) = \{x_0\}$ for some $x_0 \in X.$
Here $\phi_n$ denotes the nth iterate of $\phi.$

    The results in \cite{comp1} included not necessarily unital algebras
and not necessarily connected maximal ideal spaces. However, if we
specialize to our case (unital $B$, connected $X$), then the key elements
in the proof of Theorem 1.7 of \cite{comp1} are Lemmas 1.4 and 1.6 of
that paper. It is not hard to see that the proof of Lemma 1.4 depends
only on the fact that non-zero elements in the spectrum
of $T$ are eigenvalues of finite multiplicity, while properties (a)-(d)
of the definition of Riesz operators are all that are needed to prove
Lemma 1.6. Thus we can state the following theorem.

{\bf Theorem 1.1:} If $B$ is a unital commutative semi-simple Banach algebra
with connected maximal ideal space $X$, and if $T$ is a Riesz endomorphism
of $B$ induced by a selfmap $\phi$ of $X$, then
$\bigcap_{n=0}^\infty \phi_n(X) = \{x_0\}$  where $x_0 \in X$
is a fixed point of $\phi$ which is unique.

    If $B$ is a unital commutative semi-simple Banach algebra with maximal
ideal space $X$, for $x,y \in X$ we let
$$\|x-y\|=\sup\{|\hat{f}(x)-\hat{f}(y)|:f \in B, \|f\| \leq 1\}.$$
That is, $\|x - y\|$ is the norm of $x -y$ regarded as an element of
the dual space $B^*$ of $B$. Further, for $\varepsilon > 0$ and $a \in X$,
we let $B(a,\varepsilon)=\{x \in X: \|x - a\| < \varepsilon\}.$

    We now prove the following theorem for arbitrary unital commutative
semi-simple Banach algebras.
We remark that a related argument was used by L. Zheng
(\cite{zh} Lemma 2)  in connection with Riesz composition operators on
$H^\infty$ of the unit disk.

{\bf Theorem 1.2:} Suppose that $B$ is a unital commutative semi-simple
Banach with connected maximal ideal space $X$. Let $T$ be a Riesz
endomorphism  induced by a selfmap $\phi$ of $X$. Suppose that
$\{x_0\} = \bigcap_{n=0}^\infty \phi_n(X).$ Then for each $\varepsilon > 0$
there exists a positive integer $N$ such that $\phi_N(X) \subset
B(x_0,\varepsilon).$

{\bf Proof:} Assume that $B, X, T, \phi, x_0$ are as described in the
hypothesis. Suppose that there exists $\varepsilon > 0$ such that for
all positive integers $k$ we have that $\phi_k(X)$ is not contained in
$B(x_0,\varepsilon).$ We show that this leads to a contradiction that
$T$ is  a Riesz operator. To this end, let $x_n \in X$ satisfy
$\phi_n(x_n) \not \in B(x_0,\varepsilon),$ i.e.
$\|\phi_n(x_n) - x_0\| \geq \varepsilon.$
 From Theorem 1.1 , if $\cal U$ is a w*-open
neighborhood of $x_0$, then $\phi_n(X) \subset \cal U$ for large $n$.
In particular, $\phi_n(x_n) \in \cal U$ for large $n$, and
so it follows that $\phi(x_n) \rightarrow x_0$ in the w*-topology of $X.$

    Let $y_n=\phi_{n-1}(x_n).$ Then $y_n \rightarrow x_0$ in the
w*-topology of $X$ and $\|\phi(y_n) - x_0 \| \geq \varepsilon.$
Choose $f_n \in B$ with $\|f_n\|_B=1$ and
$|\hat{f}_n(\phi(y_n))-\hat{f}_n(x_0)| > \varepsilon - \frac{1}{n}.$
Let $K$ be any compact linear map from $B \rightarrow B.$ Then there
exist $g \in B$ and a subsequence $\{f_{n_j}\}$ with $Kf_{n_j} \rightarrow g$
in norm.

    Hence
\[\|T - K\|\geq \|Tf_{n_j} - Kf_{n_j}\| \geq \|Tf_{n_j} - g\| - \|Kf_{n_j}-g\|\]
\[\geq |\hat{f}_{n_j}(\phi(y_{n_j})) -\hat{g}(y_{n_j})| - \|Kf_{n_j} - g\|.\]

   Therefore
\[\|T - K\| \geq \limsup|\hat{f}_{n_j}(\phi(y_{n_j}))-\hat{g}(y_{n_j})|=
\limsup|\hat{f}_{n_j}(\phi(y_{n_j})) - \hat{g}(x_0)|.\]

    Also, evaluating at $x=x_0$ gives that
\[\|T-K\| \geq \limsup[|\hat{f}_{n_j}(\phi(x_0))-\hat{g}(x_0)|-\|Kf_{n_j}-g\|]=
\limsup|\hat{f}_{n_j}(x_0) - \hat{g}(x_0)|.\]

    Adding, we obtain
\[2\|T-K\| \geq \limsup[|\hat{f}_{n_j}(\phi(y_{n_j})-\hat{g}(x_0)|+
|\hat{f}_{n_j}(x_0)-\hat{g}(x_0)|],\]
and so
\[2\|T-K\| \geq \limsup|\hat{f}_{n_j}(\phi(y_{n_j}))-\hat{f}_{n_j}(x_0)|
\geq \varepsilon.\]
Therefore, $\|T - K\| \geq \varepsilon/2$ for all compact operators $K$.

   Next for each positive integer $m$, consider the endomorphism $T^m$
which is induced by $\phi_m.$ Then $T^m$ is Riesz and we have
$\dis \bigcap \phi_{mn}(X)=\{x_0\}.$
If we apply the previous argument to $T^m$, we get that $\|T^m-K\|
\geq \varepsilon/2$ for each positive integer $m$ and all compact
operators $K$. Therefore
$\dis \lim_m(\inf\{\|T^m-K\|: K$ is compact$\})^{1/m} \geq 1$.
This contradicts the assumption that $T$ is a Riesz operator. Hence we
have that if $T$ is a Riesz endomorphism induced by $\phi$, then for
each $\varepsilon > 0$ there exists a positive integer $N$ with
$\phi_N(X) \subset B(x_0,\varepsilon).$ \QED

   There are several immediate consequences of this theorem.

{\bf Corollary 1.3:} Suppose that $B$ is a unital commutative semi-simple
Banach with connected maximal ideal space $X$. Let $T$ be a Riesz
endomorphism  induced by a selfmap $\phi$ of $X$. Suppose that
$\{x_0\} = \bigcap_{n=0}^\infty \phi_n(X).$
If $x_0$ is an isolated point of $X$ in the norm topology, then
$T^Nf=\hat{f}(x_0)1$ for some positive integer $N$. \QED

{\bf Proof:} Since $x_0$ is an isolated point in the norm topology,
there exists $\varepsilon > 0$ such that $B(x_0,\varepsilon)=\{x_0\}$
and the result follows from Theorem 1.2. \QED

{\bf Corollary 1.4:} Suppose that $B$ is a unital commutative semi-simple
Banach with connected maximal ideal space $X$. Let $T$ be a Riesz
endomorphism  induced by a selfmap $\phi$ of $X$. Suppose that
$\{x_0\} = \bigcap_{n=0}^\infty \phi_n(X).$ If $B$ has no non-zero
point derivations at $x_0$, then $T^Nf=\hat{f}(x_0)1$
for some positive integer $N$. 

{\bf Proof:} Theorem 1.6.2 of \cite{br} asserts that if there are no
non-zero point derivations at $x_0,$ then $x_0$ is an isolated point of $X$
in the norm topology. The result then follows from the preceding
corollary. \QED

\begin{center}
Part 2: Dales-Davie algebras
\end{center}

    Let $X$ be a perfect compact subset of the complex plane and
let $D^\infty(X)$ denote the set of infinitely differentiable functions
on $X$. Suppose, too, that $(M_n)$ is a sequence of positive numbers
satisfying $M_0=1$ and $\dis \frac{M_{n+m}}{M_n M_m} \geq{{n+m}\choose{m}}$.
Finally, let
\[D(X, M)=\{ f \in D^{\infty}(X):\|f\|_{D}=\sum_{n=0}^{\infty}
\frac{\|f^{(n)}\|_{\infty}}{M_n} < \infty\}.\]
With pointwise addition and multiplication, $D(X,M)$ is a normed algebra.
We call such algebras Dales-Davie algebras. See \cite{dales} and \cite{dd}
for examples and basic facts about Dales-Davie algebras, and
\cite{jfhk}, \cite{jlms}
 and \cite{blaub} for some results about endomorphisms of
these algebras. We will also assume that a weight sequence $(M_n)$ is
nonanalytic meaning that $\lim_{n\to\infty} (n!/M_n)^{1/n}=0.$
Suppose that $(M_n)$ is nonanalytic,
$D(X,M)$ is a Banach algebra, and that the maximal ideal space of $D(X,M)$
is precisely $X$. In such cases, every unital endomorphism $T$
of $D(X,M)$ has the form $Tf(x)=f(\phi(x))$ for some continuous
selfmap $\phi$ of $X$. As a final definition, a selfmap $\phi$ of
a compact subset of the plane is called analytic if
$$\dis \sup_{k}\left({\frac{\|\phi^{(k)}\|_\infty}{k!}}\right)^{1/k} < \infty.$$

    Before proceeding to our result about Riesz endomorphisms of
Dales-Davie algebras, we point out the following easily verified fact.

{\bf Theorem 2.1:} If $T$ is a unital Riesz endomorphism of a unital
commutative semi-simple Banach algebra $B$ with connected maximal
ideal space $X$, then $\{f:Tf=f\}$ is a one dimensional Banach algebra.
That is, the eigenvalue $1$ has multiplicity $1$ and $Tf=f$ implies
that $f$ is a constant.

{\bf Lemma 2.2:} Let $X$ be a connected perfect compact subset of the
complex plane, $(M_n)$ a non-analytic weight sequence and $D(X,M)$
a Banach algebra with maximal ideal space $X$. Suppose that $T$ is
a Riesz endomorphism of $D(X,M)$ induced by the selfmap $\phi$ of $X$.
If $x_0$ is the fixed point of $\phi$, then $|\phi'(x_0)| < 1.$

{\bf Proof:} Assume that $D(X,M)$, $T$, $\phi$ and $x_0$ are as described.
Let $M_{x_0}=\{f:f(x_0)=0\}$ and $T_0=T_{|M_{x_0}}.$ The operator $T_0$
is a Riesz endomorphism of $M_{x_0}.$ We first show that if
$|\phi'(x_0)|=1,$ then for some positive integer $N$, $\lambda=1$
is in the spectrum $\sigma(T_0^N)$ of $T_0^N.$ Letting $\|S\|_{sp}$
denote the spectral radius of an operator $S$, we claim that if
$|\phi'(x_0)|=1,$ then $\|T_0\|_{sp}=1.$ To this end, assume that
$|\phi'(x_0)|=1$ and let $f(x)=x-x_0$. Clearly $f \in D(X,M).$ Then for each
positive integer $n$,
\[\|f\|_D \|T_0^n\| \geq \|T_0^n f\|_D \geq \|\phi_n-x_0\|_D \geq
\frac{|\phi_n'(x_0)|}{M_1}=\frac{|\phi'(x_0)|^n}{M_1}.\]
Therefore, $\dis \|T_0\|_{sp} =  \lim \|T_0^n\|^{1/n} \geq |\phi'(x_0)|=1.$
On the other hand, since the set of eigenvalues of $T_0$ is closed under
multiplication, $\sigma(T_0) \subseteq \{\lambda:|\lambda| \leq 1\}.$
Hence, $\|T_0\|_{sp}=1.$ Also, since $T_0$ is a Riesz endomorphism
every non-zero element in $\sigma(T_0)$ is an eigenvalue.
Consequently there exists
an eigenvalue $\lambda$ of $T_0$ of magnitude $1$. Again using the fact that
the eigenvalues are closed under multiplication, and the fact
that there are only finitely many eigenvalues of $T_0$ on the unit circle,
we conclude that for some positive integer $N$, $1 \in \sigma(T_0^N).$
Thus, for some non-zero $u$ in $M_{x_0}$, $T_0^N u=u.$ Further, the
function $u$ is also an eigenvector of $T^N$ on $D(X,M)$. Hence from
Theorem 2.1 $u$ is a constant function, and since $u \in M_{x_0}$,
$u$ must be $0$, a contradiction. Therefore if $T$ is a Riesz endomorphism
of the Banach algebra $D(X,M)$, then $|\phi'(x_0)| < 1.$ \QED

   Suppose that $X$ is a connected perfect compact subset of the complex
plane and $(M_n)$ is a nonanalytic weight sequence. Suppose further that
$D(X,M)$ is a Banach algebra with maximal ideal space $X$. It was
shown in \cite{jlms} that if $T$ is an endomorphism induced by an analytic
selfmap $\phi$ of $X$ and if $\|\phi'\|_\infty < 1,$ then $T$ is a compact
endomorphism. The next theorem  follows easily from this and Lemma 2.2.

{\bf Theorem 2.3:} Let $X$ be a connected perfect compact subset of the
complex plane and $(M_n)$ a nonanalytic weight sequence. If $D(X,M)$
is a Banach algebra with maximal ideal space $X$, then every Riesz
endomorphism induced by an analytic selfmap $\phi$ of $X$ is power compact.

{\bf Proof:} Suppose that $T$ is a Riesz endomorphism of $D(X,M)$ induced
by an analytic selfmap $\phi$ of $X$ and $x_0$ is the fixed point of $\phi$. From
Lemma 2.2, $|\phi'(x_0)|=c_1 < 1.$ Let $0 < c_1 < c < 1$ and let
$\cal U$ be a neighborhood of $x_0$ for which $|\phi'(t)| < c < 1$ for
all $t \in \cal U.$ Then there exists a positive integer $N$ such that
$\phi_N(x) \in \cal U$ for all $x \in X.$ Then for $n > N,$
\[|\phi_n'(x)|=|\phi'(\phi_{n-1}(x))\cdots \phi'(\phi_{N+1}(x))\phi'(\phi_N(x))
\cdots \phi'(\phi(x))\phi'(x)|.\]
Hence if $n > N$, then for all $x \in X$, $\dis |\phi_n'(x)| < c^{n-N}\|\phi'\|_\infty^N.$ Thus $\|\phi_n'\|_\infty < 1$ for large $n$, whence
$\phi$ induces  compact endomorphism of $D(X,M)$ for large $n$. Therefore
$T$ is power compact. \QED

     Further, the spectrum $\sigma(T)$ of a Riesz endomorphism $T$
of $D(X,M)$ is easy to determine in many cases. It is well known
\cite{dd} that a sufficient condition for $D(X,M)$ to be a Banach
algebra is that $X$ be uniformly regular, meaning that for all $z,w \in X$, 
there is a rectifiable arc in $X$ joining $z$ to $w$, and the metric
given by the geodesic distance between the points of $X$ is uniformly
equivalent to the Euclidean metric. Moreover, the proofs of Theorem 2.4
of \cite{jlms} when $X=[0,1]$ and Theorem 11 of \cite{edm} when $X$ is
uniformly regular show that 

\begin{center}
$\dis \sigma(T)=\{(\phi'(x_0)^n:$ n is a positive integer$\} \cup
\{0,1\}$
\end{center}
holds whenever the non-zero spectrum contains only eigenvalues. Thus
we have the following theorem.

{\bf Theorem 2.4:} Let $X$ be a uniformly regular compact subset of
the complex plane and $(M_n)$ a nonanalytic weight sequence. If $T$
is a Riesz endomorphism of $D(X,M)$ induced by $\phi$ and if $x_0$
is the fixed point of $\phi$, then
\begin{center}
$\dis \sigma(T)=\{(\phi'(x_0)^n:$ n is a positive integer$\} \cup
\{0,1\}.$
\end{center}
Also,  each nonzero element in $\sigma(T)$ has multiplicity $1.$

\begin{center}
Part 3. $C^1[0,1]$
\end{center}

   We next look at the Banach algebra $C^1[0,1]$ for examples of Riesz
endomorphisms which are not power compact. The technique in the proof
of the theorem is a slight variation of Theorem 1.2.

{\bf Theorem 3.1:} Suppose that $T$ is a unital endomorphism of $C^1[0,1]$
induced by the selfmap $\phi$ of $[0,1].$ Then $T$ is a Riesz endomorphism
if, and only if, $\dis \bigcap_{n=0}^\infty \phi_n([0,1]) =\{x_0\}$
for some $x_0 \in [0,1]$ and $\phi'(x_0) = 0.$

{\bf Proof:} First suppose that $\phi$ induces a Riesz endomorphism $T$
of $C^1[0,1]$,  $\dis \bigcap_{n=0}^\infty \phi_n([0,1]) =\{x_0\}$
and $\phi'(x_0) \neq 0.$  We aim to show that this leads to a contradiction.
To this end, choose $x_n \in [0,1]$ and  $f_n \in C^1[0,1]$
with $x_n \rightarrow x_0$,
$f_n'(\phi(x_n))=-1$, $f_n'(x_0)=1$ and $\|f_n\|_{C^1} \leq 1 + \frac{1}{n}.$
Suppose that $K$ is a compact operator on $C^1[0,1]$. Then there exist
$g \in C^1[0,1]$ and a subsequence $\{f_{n_j}\}$ such that
$Kf_{n_j} \rightarrow g.$  As in the proof of Theorem 1.2,
\[\|f_{n_j}\|\|T-K\| \geq \|Tf_{n_j} - Kf_{n_j}\| \geq \|Tf_{n_j}-g\|-
\|Kf_{n_j} - g\|,\]
and so 
\[\|T-K\| \geq \limsup_j \|Tf_{n_j}-g\|.\]
    Now
\[\limsup_j\|Tf_{n_j}-g\|_{C^1}=\]
\[\limsup_j\left(\sup_x|f_{n_j}(\phi(x))-g(x)|+\sup_x|f_{n_j}'(\phi(x))\phi'(x)
-g'(x)|\right).\]
    First evaluating at $x_{n_j}$, we get
\[\limsup_j\|Tf_{n_j} - g\|_{C^1} \geq \limsup_j[0+|-\phi'(x_{n_j})-
g'(x_{n_j})|]=\]
\[\limsup_j[|\phi'(x_{n_j})+g'(x_{n_j})|
]=|\phi'(x_0)+g'(x_0)|.\]
    Then evaluating at $x_0$, we get
\[\limsup_j \|Tf_{n_j}-g\|_{C^1} \geq \limsup_j|\phi'(x_0)-g'(x_0)|
=|\phi'(x_0)-g'(x_0)|.\]
    Adding we get $\|T - K\| \geq |\phi'(x_0)|/2.$ We recall that
the compact operator $K$ is arbitrary.

    Again, for each positive integer $m$, $T^m$ is  a Riesz
endomorphism which is induced by $\phi_m$; 
hence the preceding argument goes through for $T^m$.
That is, $\|T^m-K\| \geq |\phi'(x_0)|^m/2$ for all positive
integers $m$ and all compact operators $K$. Thus
$\lim[\inf\{\|T^m-K\|: K$ is a compact operator$\}]^{1/m} \geq |\phi'(x_0)| >0.$
This is a contradiction to the assumption that $T$ is a Riesz endomorphism.

    Conversely, assume that $\dis \bigcap_{n=0}^\infty \phi_n([0,1])=\{x_0\}$
and $\phi'(x_0)=0.$ Let $Lf=f(x_0)1$, a compact operator from $C^1[0,1]$ to 
$C^1[0,1]$. We show that $\dis \lim_n \|T^n - L\|^{1/n} = 0.$ Indeed,
let $\varepsilon > 0.$ Since  $\dis \bigcap_{n=0}^\infty \phi_n([0,1])=\{x_0\}$
and $\phi'(x_0)=0$, there exists a positive integer $N$ and a positive
number $M$ such that $|\phi_n'(t)| < M \varepsilon^n$ for $n > N$,
all $t \in [0,1].$ Then
\[|f(\phi_n(x))-f(x_0)| \leq \|f'\|_\infty\|\phi_n'\|_\infty|x-x_0| <M_1
\varepsilon^n\]
some $M_1$, all $x$, large $n$. Also
\[|f'(\phi_n(x))\phi_n'(x)| \leq \|f'\|_\infty |\phi_n'(x)| \leq M_2
\varepsilon ^n\]
for some $M_2$, all $x$, large $n$.
Therefore, $\dis \|T^n-L\| < (M_1+M_2)\varepsilon^n$, which implies
that
 $\dis \lim_n \|T^n - L\|^{1/n} = 0$ whence $T$
is a Riesz operator by Theorem 52.2 of \cite{heus}. \QED

{\bf Corollary 3.2:} There exist Riesz endomorphisms of $C^1[0,1]$ which
are not power compact.

{\bf Proof:} Every compact endomorphism of $C^1[0,1]$ has the form
$Lf=f(x_0)1$ for some $x_0 \in [0,1].$ (\cite{comp1}, Theorem 2.3).
As an example, if $\phi(x)=x^2/2$, then $\phi$ induces a Riesz endomorphism
of $C^1[0,1]$ which is not power compact. \QED
\newpage
\begin{center}
Part 4: Uniform algebras
\end{center}

   Thus far we have examples of algebras for which every Riesz endomorphism
is power compact, as well as an example of an algebra with a non-power
compact Riesz endomorphism. As we now show, every Riesz endomorphism
of the disk algebra  or $H^\infty(\Delta)$, $\Delta$ the unit disk,
is power compact. One might think that this was the case for every
uniform algebra. However, in this section we also construct an example
of a uniform algebra whose maximal ideal space is connected and which has
a non-power compact Riesz endomorphism.

We will assume some knowledge of the standard theory of uniform algebras: for
more details we refer the reader to \cite{GamBook}.

Let $A$ be a uniform algebra with maximal ideal space $X$ and let $T$ be an endomorphism of $A$ induced by a selfmap $\phi$ of $X$. It is standard that $T$ is
compact if and only if $\phi(X)$ is a norm (Gleason) compact subset of $X$
(\cite[Theorem 1]{Gamelin}).
Note, however, that the \lq if' part fails for more general Banach function algebras.

    It follows from our earlier Theorem 1.2 that if 
$\phi$ induces a Riesz endomorphism  $T$
of a uniform algebra with connected maximal ideal space $X$, 
then for some positive integer $N$,
$\phi_N(X)$ is a Gelfand (and hence also Gleason) closed subset of a closed 
norm ball of radius less than $2$.
This ball is, of course,
contained in exactly one Gleason part. 

Recall that a uniform algebra $A$ is a 
{\it unique representing measure algebra} ({\it URM-algebra}) if
every character of $A$ has a unique representing measure on the Shilov 
boundary of $A$. The disk algebra, $H^\infty(\Delta)$, and the trivial uniform algebras
$C(X)$ are all examples of URM-algebras. 
It is standard that the Gleason parts of URM-algebras are 
either one point parts or else analytic disks. 
In particular, for these 
algebras,
every closed Gleason ball of radius less than $2$ is Gleason compact.

 From the above discussion we now see immediately that
if $A$ is a URM-algebra with
connected maximal ideal space, then every Riesz endomorphism of $A$
is power-compact.

\def\N{{\bf N}}
\def\C{{\bf C}}
\def\D{{\bf D}}
See \cite{kl} for related results and examples concerning compact
endomorphisms.

In view of the above results, in order to construct non-power compact Riesz
endomorphisms of
uniform algebras, we should look at examples where the small closed Gleason
balls are
not Gleason compact. 
Klein \cite{kl} considered such uniform algebras for similar reasons.
We begin by proving a lemma which gives a sufficient condition for an
endomorphism of a uniform algebra
to be Riesz.

{\bf Lemma 4.1:}
Let $A$ be a uniform algebra on a compact space $X$ and let $x_0 \in X$.
Let $\phi$ be a self-map of $X$ which induces an endomorphism $T$ of $A$, and
as usual let $\phi_n$ be the $n$th iterate of $\phi$. Using the norm (Gleason)
distance in the maximal ideal space,
set $C_n=\sup\{\|\phi_n(x)-x_0\|: x \in X\}$.
Suppose that $C_n^{1/n} \to 0$ as $n \to \infty$. 
Then $T$ is a Riesz
endomorphism of $A$.

{\bf Proof:}
Consider the distance from $T^n$ to the compact endomorphism $L$ given by
$Lf = f(x_0) 1$. For $f \in A$ and $x \in X$ we have
$(T^n f -Lf)(x) = f(\phi_n(x)) - f(x_0)$
and it then follows from this that
$\|T^n -L\| \leq C_n$.
Thus if $C_n^{1/n} \to 0$ as $n \to \infty$, then 
$\|T^n-L\|^{1/n} \to 0$ as $n \to \infty$, and the result follows.
\QED

This condition is far from necessary, as is shown by, for example, the compact endomorphism of the
disk algebra induced by the self-map $\phi(z)=z/2$.

   We now use Lemma 4.1 to construct a non-power compact
Riesz endomorphism
of a Banach algebra of analytic functions on the unit ball of a
Banach space.

  For our example let $\cal X$ denote the closed unit ball of the complex
Banach space $E=\ell^\infty(\N^2)$, with the weak * topology. Then
 $\cal X$ is
compact and metrizable. Elements of $\cal X$ will be denoted by
$x=(x_{j,k})_{j,k=1}^\infty.$
Next define $\cal A$ as the uniform algebra
on $\cal X$ generated by the co-ordinate projections $p_{j,k}$. Then
$\cal A$ is a
subalgebra of the
uniform algebra of all continuous functions on $\cal X$ which are disk algebra
functions in each variable separately.

We consider the norm metric $d_\infty$ on the
Banach space $E$, and also the Gleason distance 
(from the norm of $\cal A^*$) on $\cal X$.
Further we  denote by $x^0$ the zero element of $E$: 
obviously $x^0 \in \cal X$.

Although we do not use this fact, it is easy to show that the maximal
ideal space of $\cal A$ is $\cal X$. Indeed, every character is determined by
what it does to all of the coordinate functionals, and must agree at all
of these with some evaluation
character at a point of $\cal X$.

 We start with a lemma to help us estimate the Gleason distance from $x^0$.

{\bf Lemma 4.2:} Let $\cal A$ and $\cal X$ be as described.

(i) If  $\alpha \in \C$ with $|\alpha|\leq 1$, 
then the selfmap $\psi$ of $\cal X$ defined by
$\psi(x) = \alpha x$ induces an endomorphism of $\cal A$.

(ii) For each $x \in \cal X$, the  Gleason distance from $x$ to $x^0$ is at
most $2 d_\infty(x,x^0)$.

{\bf Proof:} (i) For all $j,k$ we have $p_{j,k}\circ \psi = \alpha p_{j,k} \in
\cal A.$ Thus the closed subalgebra 
$\{f \in C(\cal X)$:$f \circ \psi \in \cal A\}$
of $C(\cal X)$
contains all the $p_{j,k}$ and hence all of $\cal A$, as required.

(ii) Let $x \in \cal X$ and 
assume that $x \neq x^0$. Set $R=1/d_\infty(x,x^0)$.
Let $f$ be a function in the unit ball of the
dense subalgebra of $A$ generated by $1$ and the coordinate
projections $p_{j,k}$ and
consider the function $g$ defined on the closed unit disk by
$g(z)=f(zRx)-f(x^o)$. This is simply a polynomial in $z$ which vanishes 
at $0$, and $\|g\|_\infty \leq 2 \|f\|_\infty$. Note the norm of $g$ 
is taken on the closed unit disc, so by Schwartz's Lemma
we must have $|g(z)| \leq 2\|f\|_\infty |z|$ for
$|z| \leq 1$. In particular, setting $z=1/R$ gives us
$|f(x)-f(x^0)| \leq 2 \|f\|_\infty /R= 2 d_\infty(x,x^0)$.
The result now follows.
\QED

{\bf Theorem 4.3:}
With $\cal A$ and $\cal X$ as above, there is a Riesz endomorphism of $\cal A$ 
which is
not power compact.

{\bf Proof:}
We define a certain \lq weighted shift' $\phi$
on $\cal X$ and show that this induces an endomorphism with the desired properties.

For $x\in \cal X$, we define $\phi$ by
\[
(\phi(x))_{j,k}= x_{j,k+1}/(k+1).
\]
Note that
$$
d_\infty(\phi_n(x),x^0) \leq 1/(n+1)!\eqno{(*)}
$$
 for all $n$ and all
$x \in \cal X$, and in particular $\phi$ is a selfmap of $\cal X$.

We next show that $\phi$ induces an endomorphism of $\cal A$. By definition of
$\phi$ we have,
for all $j,k$, $p_{j,k} \circ \phi = p_{j,k+1}/(k+1)$.
Thus (as for $\psi$ above) $f \circ \phi \in \cal A$  for all $f \in \cal A.$
Let $T$ be the endomorphism induced by $\phi$.
It follows easily from the preceding two lemmas and $(*)$  that $T$ is Riesz. 
To show that $T$ is not power compact, 
let $n \in \N$
and note that 
 for all $j$, $(T^n p_{j,1}) = p_{j,n+1}/(n+1)!$.
Since $(T^n p_{j,1})_{j=1}^\infty$ has no convergent subsequence,
$T^n$ is not compact for all $n$, whence the Rieesz endomorphism $T$ 
is not power compact.
\QED

This algebra was also used by Klein in \cite{kl}. 
In addition, a variety of related algebras
were considered by
Galindo, Gamelin, and Lindstr\"om \cite{GGL} in relation to weakly compact homomorphisms.
It is easy to see
from their work that weakly compact endomorphisms need not be Riesz.

Since weakly
compact endomorphisms of the
disk algebra are always compact \cite{gale}
it follows that
there are Riesz
endomorphisms of the disk algebra which are not weakly compact.
Such endomorphisms are, of course,
power compact. It is easy to see, however, that the Riesz endomorphism
constructed in Theorem 4.3 has the property that no iterate of it is 
weakly compact.

We would like to thank E. Albrecht, H. G. Dales, T. Gamelin,
R. Mortini,
M. M. Neumann and
C. J. Read 
for useful conversations which have contributed to this work.

\vspace{.3in}

{\sf  School of Mathematical Sciences

 University of Nottingham

 Nottingham NG7 2RD, England

 email: Joel.Feinstein@nottingham.ac.uk

and

 Department of Mathematics

 University of Massachusetts at Boston

 100 Morrissey Boulevard

 Boston, MA 02125-3393

 email: hkamo@cs.umb.edu

\vspace{.3in}
\vspace{.3in}

This research was supported by EPSRC grant GR/S16515/01.}

\end{document}